\documentclass[12pt]{article}
\usepackage{times}
\input{amssymb.sty}
\input{times.sty}

\def\abstr#1{{
\centering{\begin{minipage}{4.8in}\footnotesize\baselineskip=10pt
\parindent=0pt #1
\end{minipage}}\par}}

\begin{document}
\thispagestyle{empty}
\footnotesize{\noindent International Journal of Mathematics, Vol. 12, No. 2
  (2001) 159-201}

\vspace*{0.88truein}

\centerline{\bf DIMENSIONAL REDUCTION, SL(2, C)-EQUIVARIANT BUNDLES} \baselineskip=13pt 
\centerline{\bf  AND STABLE HOLOMORPHIC CHAINS}  \vspace*{0.37truein}
\centerline{\footnotesize LUIS \'ALVAREZ--C\'ONSUL} \baselineskip=12pt
\centerline{\footnotesize\it Department of Mathematics, University of
  Illinois at Urbana--Champaign,} \baselineskip=10pt 
\centerline{\footnotesize\it IL 61801, USA} \baselineskip=10pt
\centerline{\footnotesize\it E-mail: lalvarez@math.uiuc.edu} \vspace*{10pt}
\centerline{\footnotesize OSCAR GARC\'IA--PRADA} \baselineskip=12pt
\centerline{\footnotesize\it Departamento de Matem\'aticas,
  Universidad Aut\'onoma de Madrid,} \baselineskip=10pt 
\centerline{\footnotesize\it 28049 Madrid, Spain} \baselineskip=10pt
\centerline{\footnotesize\it E-mail: oscar.garcia-prada@uam.es}
\vspace*{0.225truein}

\vspace*{0.21truein}
\abstr{
In this paper we study gauge theory on SL$(2,{\mathbb C})$-equivariant
bundles over $X \times {\mathbb P}^1$, where X is a compact K\"ahler manifold,
${\mathbb P}^1$ is the complex projective line, and the action of
SL$(2,{\mathbb C})$ is trivial on $X$ and standard on ${\mathbb P}^1$. We first
classify these bundles, showing that they are in correspondence 
with objects on $X$ --- that we call holomorphic chains --- consisting
of a finite number of holomorphic bundles $E_i$ and morphisms $E_i\to
E_{i-1}$. We then prove a Hitchin--Kobayashi correspondence relating the
existence of solutions to certain natural gauge-theoretic equations
and an appropriate notion of stability for an equivariant bundle and
the corresponding chain. A central tool in this paper is a dimensional
reduction procedure which allow us to go from $X\times{\mathbb P}^1$
to $X$.}

\vspace*{0.2truein}
\abstr{
 } 

\vspace*{0.4truein}
\noindent\emph{Note}: A complete PDF file can be obtained from
\emph{http://www.worldscinet.com/journals/ijm/} \\ 
See also \emph{http://www.math.uiuc.edu/$\sim$lalvarez/papers/}

\end{document}